\documentclass[11pt]{amsart}

\theoremstyle{plain}
\newtheorem{theorem}                 {Theorem}      [section]
\newtheorem{proposition}  [theorem]  {Proposition}
\newtheorem{corollary}    [theorem]  {Corollary}

\newtheorem{lemma}        [theorem]  {Lemma}

\theoremstyle{definition}

\newtheorem{remark}       [theorem]  {Remark}

\numberwithin{equation}{section}

\DeclareMathOperator{\trace}{trace}
\DeclareMathOperator{\grad}{grad}

\DeclareMathOperator{\Id}{Id}

\DeclareMathOperator{\kernel}{Ker}

\makeindex

\begin{document}

\title[Biharmonic submanifolds with parallel mean curvature ...]
{Biharmonic submanifolds with parallel mean curvature vector field in spheres}

\author{A. Balmu\c s}
\author{C. Oniciuc}
\address{Faculty of Mathematics, ``Al.I.~Cuza'' University of Iasi\\
\newline
Bd. Carol I Nr. 11 \\
700506 Iasi, ROMANIA}

\email{adina.balmus@uaic.ro, oniciucc@uaic.ro}


\thanks{The first author was supported by Grant POSDRU/89/1.5/S/49944, Romania.
}

\begin{abstract}
We present some results on the boundedness of the mean curvature of
proper biharmonic submanifolds in spheres. A partial classification
result for proper biharmonic submanifolds with parallel mean
curvature vector field in spheres is obtained. Then, we completely classify the proper biharmonic submanifolds in spheres with parallel mean curvature vector field and parallel Weingarten operator associated to the mean curvature vector field.
\end{abstract}


\maketitle

\section{Introduction}

Biharmonic maps between two Riemannian manifolds $(M,g)$ and
$(N,h)$, $M$ compact, generalize harmonic maps (see \cite{ES64}) and
represent the critical points of the {\it bienergy functional}
$$
E_2:C^\infty(M,N)\to\mathbb{R},\qquad E_{2}\left( \phi \right) =
\frac{1}{2} \int_{M} |\tau(\phi)|^{2}\, v_{g},
$$
where $\tau(\phi)=\trace\nabla d\phi$ denotes the {\it tension
field} associated to the map $\phi$. We recall that harmonic maps
are characterized by the vanishing of the tension field (see, for
example, \cite{EL83}).

The first variation of $E_2$, obtained by G.Y.~Jiang in \cite{J86},
shows that $\phi$ is a biharmonic map if and only if its {\it
bitension field} vanishes
\begin{eqnarray}\label{eq: bih_eq}
\tau_2(\phi)&=&-J(\tau(\phi))=- \Delta \tau(\phi) - \trace{R^{N}}( d\phi\cdot,
\tau(\phi) ) d\phi\cdot\nonumber\\&=&0,
\end{eqnarray}
i.e. $\tau(\phi)\in\kernel J$, where $J$ is, formally, the Jacobi
operator associated to $\phi$. Here $\Delta$ denotes the rough
Laplacian on sections of the pull-back bundle $\phi^{-1}(TN)$ and
$R^N$ denotes the curvature operator on $(N,h)$, and we use the
following sign conventions
$$
\Delta V=-\trace\nabla^2 V,\qquad \forall V\in C(\phi^{-1}(TN)),
$$
$$
R^N(X,Y)=[\nabla_X,\nabla_Y]-\nabla_{[X,Y]},\qquad \forall X,Y\in C(TN).
$$
When $M$ is not compact a map $\phi:(M,g)\to(N,h)$ is said to be
biharmonic if it is a solution of equation \eqref{eq: bih_eq}. As
$J$ is a linear operator, any harmonic map is biharmonic. We call
{\it proper biharmonic} the non-harmonic biharmonic maps, and the
submanifolds with non-harmonic (non-minimal) biharmonic inclusion
map are called {\it proper biharmonic submanifolds}.

One can easily construct proper biharmonic maps between Euclidean
spaces, for example by choosing third order polynomial maps or by
using the Almansi property (see \cite{A98}). Regarding proper
biharmonic Riemannian immersions into the Euclidean space, they are
characterized by the equation $\Delta H=0$, where $H$ denotes the
mean curvature vector field, i.e. they are also biharmonic in the
sense of Chen (see \cite{C84}).

A nonexistence result for proper biharmonic maps was obtained by
requesting a compact domain and a non-positively curved codomain
\cite{J86}. Moreover, the nonexistence of proper biharmonic
Riemannian immersions with constant mean curvature in non-positively
curved spaces was proved (see \cite{O02}). Other nonexistence
results, mainly regarding proper biharmonic Riemannian immersions
into non-positively curved manifolds can be found in
\cite{BMO10,BMO08,CMO02,BI98,D92,HV95, NU}. Surprisingly, in \cite{OT10}
the author constructed examples of proper biharmonic Riemannian
immersions (of non-constant mean curvature) in conformally flat
negatively curved spaces.

On the other hand there are many examples of proper biharmonic
submanifolds in positively curved spaces.

In this paper we study proper biharmonic submanifolds
in Euclidean spheres with additional extrinsic properties: parallel mean curvature vector field or parallel Weingarten operator associated to the mean curvature vector field, obtaining some rigidity results.

The paper is organized as follows. In the preliminary section we gather some known results on proper biharmonic submanifolds in the unit Euclidean sphere $\mathbb{S}^n$. This section also recalls the Moore decomposition lemma.

In the main section we first prove, for compact proper biharmonic
submanifolds of $\mathbb{S}^n$, a boundedness condition involving the mean curvature $|H|$ and the norm $|A_H|$ of the Weingarten operator associated to the mean curvature vector field (Theorem
\ref{th: classif_bih compact}).

Then, the proper biharmonic submanifolds with parallel mean curvature vector field in unit Euclidean spheres are studied. It is known that a constant mean
curvature proper biharmonic submanifold in $\mathbb{S}^n$ satisfies
$|H|\in(0,1]$, and $|H|=1$ if and only if it is minimal in a small
hypersphere $\mathbb{S}^{n-1}(1/\sqrt{2})$ (see \cite{O03}). We
prove here that the mean curvature of the proper biharmonic
submanifolds $M^m$ with parallel mean curvature vector field in
$\mathbb{S}^n$ takes values in $(0,\frac{m-2}{m}]\cup\{1\}$, and we
determine the proper biharmonic submanifolds with parallel mean
curvature and $|H|=\frac{m-2}{m}$ (Theorem \ref{th: classif_parallelH}).

Finally, we investigate proper biharmonic submanifolds in
spheres with parallel mean curvature vector field, parallel
Weingarten operator associated to the mean curvature vector field,
and $|H|\in(0,\frac{m-2}{m})$. We first prove that such submanifolds have exactly two distinct principal curvatures in the direction of
$H$ (Corollary \ref{cor: bih_2curv}) and then, using the Moore
Lemma, we determine all of them (Theorem \ref{th:
classif_parallel_H_AH}).

We shall work in the $C^\infty$ category, i.e. all manifolds,
metrics, connections, maps, sections are assumed to be smooth. All
manifolds are assumed to be connected.

\section{Preliminaries}

The biharmonic equation \eqref{eq: bih_eq} for the inclusion map
$\imath:M^m\to\mathbb{S}^n$ of a submanifold $M$ in $\mathbb{S}^n$ writes
$$
\Delta H=mH,
$$
where $H$ denotes the mean curvature vector field of $M$ in
$\mathbb{S}^n$. Although simple, this equation is not used in order
to obtain examples and classification results. The following
characterization, obtained by splitting the bitension field in its
normal and tangent components, proved to be more useful.

\begin{theorem}[\cite{O02}]\label{th: bih subm S^n}
{\rm (i)} The canonical inclusion $\imath:M^m\to\mathbb{S}^n$ of a
submanifold $M$ in an $n$-dimensional unit Euclidean sphere is
biharmonic if and only if
\begin{equation}\label{caract_bih_spheres}
\left\{
\begin{array}{l}
\ \Delta^\perp H+\trace B(A_H(\cdot),\cdot)-mH=0
\\ \mbox{} \\
\ 4\trace A_{\nabla^\perp_{(\cdot)}H}(\cdot)+m\grad(\vert H
\vert^2)=0,
\end{array}
\right.
\end{equation}
where $A$ denotes the Weingarten operator, $B$ the second
fundamental form, $H$ the mean curvature vector field,
$\nabla^\perp$ and $\Delta^\perp$ the connection and the Laplacian
in the normal bundle of $M$ in $\mathbb{S}^n$, and $\grad$ denotes
the gradient on $M$.

{\rm (ii)}If $M$ is a submanifold with parallel mean curvature
vector field, i.e. $\nabla^\perp H=0$, in $\mathbb{S}^n$, then $M$
is biharmonic if and only if
\begin{equation}\label{eq: caract_bih_Hparallel_I}
\trace B(A_H(\cdot),\cdot)=mH,
\end{equation}
or equivalently,
\begin{equation}\label{eq: caract_bih_Hparallel_II}
\left\{
\begin{array}{ll}
|A_H|^2=m|H|^2
\\ \mbox{} \\
\langle A_H, A_\eta\rangle=0,\quad \forall\eta\in C(NM), \eta\perp H,
\end{array}
\right.
\end{equation}
where $NM$ denotes the normal bundle of $M$ in $\mathbb{S}^n$.
\end{theorem}

We recall that the small hypersphere
\begin{equation}\label{eq: small_hypersphere}
\mathbb{S}^{n-1}(1/\sqrt 2)=\left\{(x,1/\sqrt 2)\in\mathbb{R}^{n+1}: x\in \mathbb{R}^n, |x|^2=1/2\right\}\subset\mathbb{S}^{n}
\end{equation}
and the standard products of spheres $\mathbb{S}^{n_1}(1/\sqrt 2)\times\mathbb{S}^{n_2}(1/\sqrt 2)$, given by
\begin{equation}\label{eq: product_spheres}
\left\{(x,y)\in\mathbb{R}^{n+1}: x\in\mathbb{R}^{n_1+1}, y\in\mathbb{R}^{n_2+1}, |x|^2=|y|^2=1/2\right\}\subset\mathbb{S}^{n},
\end{equation}
$n_1+n_2=n-1$ and $n_1\neq n_2$, are the main examples of proper biharmonic submanifolds in $\mathbb{S}^n$ (see \cite{CMO01, J86}). Inspired by these examples, by using their minimal submanifolds, two methods of construction for proper biharmonic submanifolds in spheres were given.


\begin{theorem}[\cite{CMO02}]\label{th: rm_minim}
Let $M$ be a submanifold in a small hypersphere $\mathbb{S}^{n-1}(1/\sqrt 2)\subset\mathbb{S}^n$. Then $M$ is proper biharmonic in $\mathbb{S}^{n}$ if and only if $M$ is minimal in
$\mathbb{S}^{n-1}(1/\sqrt 2)$.
\end{theorem}

We note that the proper biharmonic submanifolds of $\mathbb{S}^n$
obtained from minimal submanifolds of the proper biharmonic
hypersphere $\mathbb{S}^{n-1}(1/\sqrt 2)$ are pseudo-umbilical, i.e.
$A_H=|H|^2\Id$, have parallel mean curvature vector field and mean
curvature $|H|=1$. Clearly, $\nabla A_H=0$.

\begin{theorem}[\cite{CMO02}]\label{th:hipertor}
Let $n_1, n_2$ be two positive integers such that $n_1+n_2=n-1$, and
let $M_1$ be a submanifold in $\mathbb{S}^{n_1}(1/\sqrt 2)$ of
dimension $m_1$, with $0 \leq m_1 \leq n_1$, and let $M_2$ be a
submanifold in $\mathbb{S}^{n_2}(1/\sqrt 2)$ of dimension $m_2$,
with $0 \leq m_2 \leq n_2$. Then $M_1\times M_2$ is proper
biharmonic in $\mathbb{S}^n$ if and only if
\begin{equation*}
\left\{
\begin{array}{ll}
m_1\neq m_2\\
\tau_2(\imath_1)+2(m_2-m_1)\tau(\imath_1)=0\\
\tau_2(\imath_2)-2(m_2-m_1)\tau(\imath_2)=0\\
|\tau(\imath_1)|=|\tau(\imath_2)|,
\end{array}
\right.
\end{equation*}
where $\imath_1:M_1\to \mathbb{S}^{n_1}(1/\sqrt 2)$ and
$\imath_2:M_2\to \mathbb{S}^{n_2}(1/\sqrt 2)$ are the canonical inclusions.
\end{theorem}
Obviously, if $M_2$ is minimal in $\mathbb{S}^{n_2}(1/\sqrt 2)$,
then $M_1\times M_2$ is biharmonic in $\mathbb{S}^n$ if and only if
$M_1$ is minimal in $\mathbb{S}^{n_1}(1/\sqrt 2)$. The proper
biharmonic submanifolds obtained in this way are no longer
pseudo-umbilical, but still have parallel mean curvature vector
field and their mean curvature is $|H|=\frac{|m_1-m_2|}{m}\in(0,1)$,
where $m=m_1+m_2$. Moreover, $\nabla A_H=0$ and the principal
curvatures in the direction of $H$, i.e. the eigenvalues of $A_H$,
are constant on $M$ and given by
$\lambda_1=\ldots=\lambda_{m_1}=\frac{m_1-m_2}{m}$,
$\lambda_{m_1+1}=\ldots=\lambda_{m_1+m_2}=-\frac{m_1-m_2}{m}$.

In the proof of the main results of this paper we shall also use the
following lemma.
\begin{lemma}[Moore Lemma, \cite{M71}]\label{lem: Moore lemma}
Suppose that $M_1$ and $M_2$ are connected Riemannian manifolds and
that
$$
\varphi: M_1\times M_2\to \mathbb{R}^r
$$
is an isometric immersion of the Riemannian product. If the second
fundamental form $\widetilde{B}$ of $\varphi$ has the property
$$
\widetilde{B}(X,Y)=0,
$$
for all $X$ tangent to $M_1$, $Y$ tangent to $M_2$, then $\varphi$
is a product immersion
$\varphi=\varphi_0\times\varphi_1\times\varphi_2$, where
$\varphi_0:M_1\times M_2\to \mathbb{R}^{n_0}$ is constant,
$\varphi_i:M_i\to \mathbb{R}^{n_i}$, $i=1,2$, and
$\mathbb{R}^r=\mathbb{R}^{n_0}\oplus\mathbb{R}^{n_1}\oplus\mathbb{R}^{n_2}$
is an orthogonal decomposition. Moreover, $\mathbb{R}^{n_1}$ is the
subspace of $\mathbb{R}^r$ generated by all vectors tangent to
$M_1\times\{p_2\}$, for all $p_2\in M_2$, and $\mathbb{R}^{n_2}$ is
the subspace generated by all vectors tangent to $\{p_1\}\times
M_2$, for all $p_1\in M_1$.
\end{lemma}

\section{Main results}

\subsection{Compact proper biharmonic submanifolds in spheres}\quad

The following result for proper biharmonic constant mean curvature submanifolds
in spheres  was obtained.
\begin{theorem}[\cite{O03}]\label{th: classif_bih const mean}
Let $M$ be a proper biharmonic submanifold with constant mean
curvature in $\mathbb{S}^n$. Then $|H|\in(0,1]$. Moreover, if
$|H|=1$, then $M$ is a minimal submanifold of a small hypersphere
$\mathbb{S}^{n-1}(1/\sqrt{2})\subset\mathbb{S}^n$.
\end{theorem}

If the condition on the mean curvature to be constant is replaced by
the condition on the submanifold to be compact, we obtain the
following.

\begin{theorem}\label{th: classif_bih compact}
Let $M$ be a compact proper biharmonic submanifold of $\mathbb{S}^n$.
Then either
\begin{itemize}
\item[(i)] there exists a point $p\in M$ such that $|A_H(p)|^2<m|H(p)|^2$,
\end{itemize}
or
\begin{itemize}
\item[(ii)] $|A_H|^2=m|H|^2$. In this case, $M$ has parallel mean
curvature vector field and $|H|\in (0,1]$.
\end{itemize}
\end{theorem}

\begin{proof}
Let $M$ be a proper biharmonic submanifold of
$\mathbb{S}^n$. The first equation of \eqref{caract_bih_spheres} implies
that
$$
\langle \Delta^\perp H,H\rangle=m|H|^2-|A_H|^2,
$$
and by using the Weitzenb\" ock formula,
$$
\frac{1}{2}\Delta|H|^2=\langle \Delta^\perp
H,H\rangle-|\nabla^\perp H|^2,
$$
we obtain
\begin{equation}\label{eq: expr_bih}
\frac{1}{2}\Delta|H|^2=m|H|^2-|A_H|^2-|\nabla^\perp H|^2.
\end{equation}
As $M$ is compact, by integrating equation \eqref{eq: expr_bih} on
$M$ we get
\begin{equation*}\label{eq: cons_int_Weitz}
\int_M (m|H|^2-|A_H|^2)\,v_g\geq 0,
\end{equation*}
and (i) and the first part of (ii) follow. Then, it is easy to see
that
\begin{equation*}\label{eq: Cauchy A_H}
m|H|^4\leq |A_H|^2,
\end{equation*}
for any submanifold of a given Riemannian manifold, so when
$|A_H|^2=m|H|^2$ we get $|H|\in(0,1]$. Moreover, by integrating
\eqref{eq: expr_bih}, we obtain $\nabla^\perp H=0$ and we conclude
the proof.
\end{proof}

Regarding the mean curvature, from Theorem \ref{th: classif_bih
compact}, we get the following result.
\begin{corollary}
Let $M$ be a compact proper biharmonic submanifold of
$\mathbb{S}^n$. Then either
\begin{itemize}
\item[(i)] there exists a point $p\in M$ such that $|H(p)|<1$,
\end{itemize}
or
\begin{itemize}
\item[(ii)] $|H|=1$. In this case $M$ is a minimal submanifold of a small hypersphere $\mathbb{S}^{n-1}(1/\sqrt{2})\subset\mathbb{S}^n$.
\end{itemize}
\end{corollary}

\subsection{Biharmonic submanifolds with $\nabla^\perp H=0$ in spheres}\quad

The following result concerning proper biharmonic surfaces with
parallel mean curvature vector field was proved.

\begin{theorem}[\cite{BMO08}]\label{th: classif_surf_parallelH}
Let $M^2$ be a proper biharmonic surface with parallel mean
curvature vector field in $\mathbb{S}^n$. Then $M$ is minimal in a
hypersphere $\mathbb{S}^{n-1}(1/\sqrt 2)$ in $\mathbb{S}^n$.
\end{theorem}

We shall further see that, when $m>2$, the situation is more complex
and, apart from $1$, the mean curvature can assume other lower
values, as expected in view of Theorem \ref{th:hipertor}.

First, let us prove an auxiliary result, concerning non-full proper
biharmonic submanifolds of $\mathbb{S}^n$, which generalizes Theorem
5.4 in \cite{BMO08}.

\begin{proposition}\label{prop: non-full_bih_spheres}
Let $M^m$ be a submanifold of a small hypersphere
$\mathbb{S}^{n-1}(a)$ in $\mathbb{S}^n$, $a\in(0, 1)$. Then $M$ is
proper biharmonic in $\mathbb{S}^n$ if and only if either $a=1/\sqrt
2$ and $M$ is minimal in $\mathbb{S}^{n-1}(1/\sqrt 2)$, or $a >
1/\sqrt 2$ and $M$ is minimal in a small hypersphere
$\mathbb{S}^{n-2}(1/\sqrt 2)$ in $\mathbb{S}^{n-1}(a)$. In both
cases, $|H|=1$.
\end{proposition}

\begin{proof}
The converse follows immediately by using Theorem \ref{th: rm_minim}.

In order to prove the other implication, denote by $\mathbf{j}$ and
$\mathbf{i}$ the inclusion maps of $M$ in $\mathbb{S}^{n-1}(a)$ and
of $\mathbb{S}^{n-1}(a)$ in $\mathbb{S}^{n}$, respectively.

Up to an isometry of $\mathbb{S}^n$, we can consider
$$
\mathbb{S}^{n-1}(a)=\left\{(x^1,\ldots,
x^n,\sqrt{1-a^2})\in\mathbb{R}^{n+1}:
\sum_{i=1}^{n}(x^i)^2=a^2\right\}\subset\mathbb{S}^{n}.
$$
Then
$$
C(T\mathbb{S}^{n-1}(a))=\left\{(X^1,\ldots, X^{n},0)\in
C(T\mathbb{R}^{n+1}): \sum_{i=1}^{n}x^iX^i=0\right\},
$$
while $\displaystyle{\eta=\frac{1}{c}\left(x^1,\ldots,
x^{n},-\frac{a^2}{\sqrt{1-a^2}}\right)}$ is a unit section in the
normal bundle of $\mathbb{S}^{n-1}(a)$ in $\mathbb{S}^{n}$, where
$c^2=\frac{a^2}{1-a^2}$, $c>0$. The tension and bitension fields of
the inclusion $\imath=\mathbf{i}\circ \mathbf{j}:M\to\mathbb{S}^n$,
are given by
$$
\tau(\imath)=\tau(\mathbf{j})-\frac{m}{c}\eta,\qquad
\tau_2(\imath)=\tau_2(\mathbf{j})-\frac{2m}{c^2}\tau(\mathbf{j})+
\frac{1}{c}\left\{|\tau(\mathbf{j})|^2-\frac{m^2}{c^2}(c^2-1)\right\}\eta.
$$
Since $M$ is biharmonic in $\mathbb{S}^{n}$, we obtain
\begin{equation}\label{eq: tau2_j}
\tau_2(\mathbf{j})=\frac{2m}{c^2}\tau(\mathbf{j})
\end{equation}
and
$$
|\tau(\mathbf{j})|^2=\frac{m^2}{c^2}(c^2-1)=\frac{m^2}{a^2}(2a^2-1).
$$
From here $a\geq 1/\sqrt 2$.
Also,
$$
|\tau(\imath)|^2=|\tau(\mathbf{j})|^2+\frac{m^2}{c^2}=m^2.
$$
This implies that the mean curvature of $M$ in $\mathbb{S}^{n}$ is $1$.

The case $a=1/\sqrt 2$ is solved by Theorem \ref{th: rm_minim}.

Consider $a>1/\sqrt 2$, thus $\tau(\mathbf{j})\neq 0$. As $|H|=1$,
by applying Theorem \ref{th: classif_bih const mean}, $M$ is a
minimal submanifold of a small hypersphere $\mathbb{S}^{n-1}(1/\sqrt
2)\subset\mathbb{S}^{n}$, so it is pseudo-umbilical and with
parallel mean curvature vector field in $\mathbb{S}^{n}$
(\cite{C73}). From here it can be proved that $M$ is also
pseudo-umbilical and with parallel mean curvature vector field in
$\mathbb{S}^{n-1}(a)$. As $M$ is not minimal in
$\mathbb{S}^{n-1}(a)$, it follows that $M$ is a minimal submanifold
of a small hypersphere $\mathbb{S}^{n-2}(b)$ in
$\mathbb{S}^{n-1}(a)$. By a straightforward computation, equation
\eqref{eq: tau2_j} implies $b=1/\sqrt 2$ and the proof is completed.

\end{proof}

Since every small sphere $\mathbb{S}^{n'}(a)$ in $\mathbb{S}^n$, $a\in(0,1)$, is contained into a great sphere $\mathbb{S}^{n'+1}$ of $\mathbb{S}^n$, from Proposition \ref{prop: non-full_bih_spheres} we have the following.

\begin{corollary}\label{cor: non-full}
Let $M^m$ be a submanifold of a small sphere $\mathbb{S}^{n'}(a)$ in $\mathbb{S}^n$, $a\in(0, 1)$. Then $M$ is proper biharmonic in
$\mathbb{S}^n$ if and only if either $a=1/\sqrt 2$ and $M$ is
minimal in $\mathbb{S}^{n'}(1/\sqrt 2)$, or $a>1/\sqrt 2$ and $M$ is minimal in a small hypersphere $\mathbb{S}^{n'-1}(1/\sqrt 2)$ in
$\mathbb{S}^{n'}(a)$. In both cases, $|H|=1$.
\end{corollary}

Let $M^m$ be a submanifold in $\mathbb{S}^n$. For our purpose it is
convenient to define, following \cite{AdC94} and \cite{AGM10}, the
$(1,1)$-tensor field $\Phi=A_H-|H|^2 I$, where $I$ is the identity
on $C(TM)$. We notice that $\Phi$ is symmetric, $\trace \Phi=0$ and
\begin{equation}\label{eq: normPhi}
|\Phi|^2=|A_H|^2-m|H^4|.
\end{equation}
Moreover, $\Phi=0$ if and only if $M$ is pseudo-umbilical.

By using the Gauss equation of $M$ in $\mathbb{S}^n$, one gets the
curvature tensor field of $M$ in terms of $\Phi$ as follows.
\begin{lemma}\label{lem: curv_Phi}
Let $M^m$ be a submanifold in $\mathbb{S}^n$ with nowhere zero mean curvature
vector field. Then the curvature tensor field of $M$ is given by
\begin{eqnarray}\label{eq: curv_Phi}
R(X,Y)Z&=&(1+|H|^2)(\langle Z,Y\rangle X-\langle Z,X\rangle Y)\nonumber\\
&&+\frac{1}{|H|^2}(\langle Z,\Phi(Y)\rangle \Phi(X)-\langle Z,\Phi(X)\rangle \Phi(Y))\nonumber\\
&&+\langle Z,\Phi(Y)\rangle X-\langle Z,\Phi(X)\rangle Y+\langle Z,Y\rangle \Phi(X)-\langle Z,X\rangle \Phi(Y)\nonumber\\
&&+\sum_{a=1}^{k-1}\{\langle Z,A_{\eta_a}(Y)\rangle A_{\eta_a}(X)
-\langle Z,A_{\eta_a}(X)\rangle A_{\eta_a}(Y)\},
\end{eqnarray}
for all $X, Y, Z\in C(TM)$, where $\{H/|H|,\eta_{a}\}_{a=1}^{k-1}$,
$k=n-m$, denotes a local orthonormal frame field in the normal
bundle of $M$ in $\mathbb{S}^n$.
\end{lemma}

In the case of hypersurfaces, i.e. $k=1$, the previous result holds by making the convention that $\displaystyle{\sum_{a=1}^{k-1}\{\ldots\}=0}$.

For what concerns the expression of $\trace\nabla^2\Phi$, which will be needed further, the following result holds.
\begin{lemma}\label{lem:nabla_2_Phi}
Let $M^m$ be a submanifold in $\mathbb{S}^n$ with nowhere zero mean
curvature vector field. If $\nabla^\perp H=0$, then $\nabla\Phi$ is
symmetric and
\begin{eqnarray}\label{eq: trace_nabla2_Phi}
(\trace\nabla^2\Phi)(X)&=&-|\Phi|^2X+\left(m+m|H|^2-\frac{|\Phi|^2}{|H|^2}\right)\Phi(X)+m\Phi^2(X)\nonumber\\
&&-\sum_{a=1}^{k-1}\langle\Phi,A_{\eta_a}\rangle A_{\eta_a}(X).
\end{eqnarray}
\end{lemma}

\begin{proof}
From the Codazzi equation, as $\nabla^\perp H=0,$ we get
$(\nabla A_H)(X,Y)=(\nabla A_H)(Y,X)$, for all $X, Y\in C(TM)$, where
$$
(\nabla A_H)(X,Y)=(\nabla_X A_H)(Y)=\nabla_X A_H(Y)-A_H(\nabla_X Y).
$$
As the mean curvature of $M$ is constant we have $\nabla\Phi=\nabla
A_H$, thus $\nabla\Phi$ is symmetric.

We recall the Ricci commutation formula
\begin{equation}\label{eq: Ricci_comm_Phi}
(\nabla^2\Phi)(X,Y,Z)-(\nabla^2\Phi)(Y,X,Z)=R(X,Y)\Phi(Z)-\Phi(R(X,Y)Z),
\end{equation}
 for all $X,Y,Z\in C(TM)$, where
\begin{eqnarray*}
(\nabla^2\Phi)(X,Y,Z)&=&(\nabla_X\nabla\Phi)(Y,Z)\\
&=&\nabla_X((\nabla\Phi)(Y,Z))
-(\nabla\Phi)(\nabla_XY,Z)-(\nabla\Phi)(Y, \nabla_X Z).
\end{eqnarray*}

Consider $\{X_i\}_{i=1}^m$ to be a local orthonormal frame field on
$M$ and $\{H/|H|,\eta_{a}\}_{a=1}^{k-1}$, $k=n-m$, a local
orthonormal frame field in the normal bundle of $M$ in
$\mathbb{S}^n$. As $\eta_a$ is orthogonal to $H$, we get $\trace
A_{\eta_a}=0$, for all $a=1,\ldots,k-1$. Using also the symmetry of
$\Phi$ and $\nabla\Phi$, \eqref{eq: Ricci_comm_Phi} and \eqref{eq:
curv_Phi}, we have
\begin{eqnarray*}
(\trace\nabla^2\Phi)(X)&=&\sum_{i=1}^m(\nabla^2\Phi)(X_i,X_i,X)=\sum_{i=1}^m(\nabla^2\Phi)(X_i,X,X_i)\\
&=&\sum_{i=1}^m\{(\nabla^2\Phi)(X,X_i,X_i)+R(X_i,X)\Phi(X_i)-\Phi(R(X_i,X)X_i)\}\\
&=&\sum_{i=1}^m(\nabla^2\Phi)(X,X_i,X_i)\\
&&-|\Phi|^2X+\left(m+m|H|^2-\frac{|\Phi|^2}{|H|^2}\right)\Phi(X)+m\Phi^2(X)\nonumber\\
&&+\sum_{a=1}^{k-1}\{(A_{\eta_a}\circ\Phi-\Phi\circ A_{\eta_a})(A_{\eta_a}(X))
-\langle\Phi,A_{\eta_a}\rangle A_{\eta_a}(X)\}.
\end{eqnarray*}
By a straightforward computation,
$$
\sum_{i=1}^m(\nabla^2\Phi)(X,X_i,X_i)=\nabla_X(\trace\nabla\Phi)=\nabla_X\grad(\trace\Phi)=0.
$$
Moreover, from the Ricci equation, since $\nabla^\perp H=0$,
we obtain $A_{\eta_a}\circ A_H=A_H\circ A_{\eta_a}$, thus
$A_{\eta_a}\circ\Phi=\Phi\circ A_{\eta_a}$, and we end
the proof of this lemma.

\end{proof}

We shall also use the following lemma.
\begin{lemma}\label{lem: cons_Weitz}
Let $M^m$ be a submanifold in $\mathbb{S}^n$ with nowhere zero mean curvature
vector field. If $\nabla^\perp H=0$ and $A_H$ is orthogonal to $A_{\eta_a}$,
for all $a=1,\ldots,k-1$, then
\begin{equation}\label{eq: cons_Weitz}
-\frac{1}{2}\Delta|\Phi|^2=|\nabla\Phi|^2+\left(m+m|H|^2
-\frac{|\Phi|^2}{|H|^2}\right)|\Phi|^2+m\trace\Phi^3.
\end{equation}
\end{lemma}

\begin{proof}
Since $A_H$ is orthogonal to $A_{\eta_a}$ and $\trace A_{\eta_a}=0$,
we get $\langle\Phi,A_{\eta_a}\rangle=0$, for all $a=1,\ldots,k-1$,
and \eqref{eq: trace_nabla2_Phi} becomes
\begin{equation}\label{eq: trace_nabla2_Phi_II}
(\trace\nabla^2\Phi)(X)=-|\Phi|^2X+\left(m+m|H|^2-\frac{|\Phi|^2}{|H|^2}\right)\Phi(X)+m\Phi^2(X).
\end{equation}

Now, the Weitzenb\"{o}ck formula,
$$
-\frac{1}{2}\Delta|\Phi|^2=|\nabla\Phi|^2+\langle \Phi, \trace\nabla^2\Phi\rangle,
$$
together with the symmetry of $\Phi$ and \eqref{eq:
trace_nabla2_Phi_II}, leads to the conclusion.
\end{proof}

We also recall here the Okumura Lemma.
\begin{lemma}[Okumura Lemma, \cite{O74}]\label{lem: Okumura}
Let $b_1,\ldots, b_m$ be real numbers such that
$\displaystyle{\sum_{i=1}^m b_i=0}$. Then
$$
-\frac{m-2}{\sqrt{m(m-1)}}\left(\sum_{i=1}^m b_i^2\right)^{3/2}\leq
\sum_{i=1}^m b_i^3\leq\frac{m-2}{\sqrt{m(m-1)}}\left(\sum_{i=1}^m b_i^2\right)^{3/2}.
$$
Moreover, equality holds in the right-hand (respectively, left-hand)
side if and only if $(m - 1)$ of the $b_i$'s are nonpositive
(respectively, nonnegative) and equal.
\end{lemma}

By using the above lemmas we obtain the following result on the
boundedness of the mean curvature of proper biharmonic submanifolds
with parallel mean curvature in spheres, as well as a partial
classification result. We shall see that $|H|$ does not fill out all
the interval $(0,1]$.

\begin{theorem}\label{th: classif_parallelH}
Let $M^m$, $m>2$, be a proper biharmonic submanifold with
parallel mean curvature vector field in $\mathbb{S}^n$ and $|H|\in(0,1)$. Then $|H|\in(0,\frac{m-2}{m}]$. Moreover, $|H|=\frac{m-2}{m}$ if and only if $M$ is an open part of
a standard product
$$
M_1^{m-1}\times \mathbb{S}^{1}(1/\sqrt 2)\subset\mathbb{S}^{n},
$$
where $M_1$ is a minimal submanifold in $\mathbb{S}^{n-2}(1/\sqrt 2)$.
\end{theorem}

\begin{proof}
Consider the tensor field $\Phi$ associated to $M$. Since it is
traceless, Lemma \ref{lem: Okumura} implies that
\begin{equation}\label{eq: cons_Okumura}
\trace\Phi^3\geq-\frac{m-2}{\sqrt{m(m-1)}}|\Phi|^3.
\end{equation}
By \eqref{eq: caract_bih_Hparallel_II}, as $M$ is proper biharmonic
with parallel mean curvature vector field, $|A_H|^2=m|H^2|$ and
$\langle A_H, A_\eta\rangle=0$, for all $\eta\in C(NM)$, $\eta$
orthogonal to $H$. From \eqref{eq: normPhi} we obtain
\begin{equation}\label{eq: normPhi_bih}
|\Phi|^2=m|H|^2(1-|H|^2),
\end{equation}
thus $|\Phi|$ is constant. We can apply Lemma \ref{lem: cons_Weitz}
and, using \eqref{eq: cons_Okumura} and \eqref{eq: normPhi_bih},
equation \eqref{eq: cons_Weitz} leads to
$$
0\geq m^2|H|^3(1-|H|^2)\left(2|H|-\frac{m-2}{\sqrt{m-1}}\sqrt{1-|H|^2}\right),
$$
thus $|H|\in (0,\frac{m-2}{m}]$.

The condition $|H|=\frac{m-2}{m}$ holds if and only if
$\nabla\Phi=0$ and we have equality in \eqref{eq: cons_Okumura}.
This is equivalent to the fact that $\nabla A_H=0$ and, by the
Okumura Lemma, the principal curvatures in the direction of $H$ are
constant functions on $M$ and given by
\begin{align}\label{eq: curv_Okumura}
\lambda_1=\ldots=\lambda_{m-1}&=\lambda=\frac{m-2}{m},\nonumber\\
\lambda_m&=\mu=-\frac{m-2}{m}.
\end{align}
Further, we consider the distributions
\begin{align*}
T_\lambda&=\{X\in TM: A_H(X)=\lambda X\},\qquad \dim T_\lambda=m-1,\\
T_\mu&=\{X\in TM: A_H(X)=\mu X\},\qquad \dim T_\mu=1.
\end{align*}
One can easily verify that, as $A_H$ is parallel, $T_\lambda$ and
$T_\mu$ are mutually orthogonal, smooth, involutive and parallel,
and the de Rham decomposition theorem (see \cite{KN63}) can be
applied. 

Thus, for every $p_0\in M$ there exists a neighborhood
$U\subset M$ which is isometric to a product
$\widetilde{M}_1^{m-1}\times I$, $I=(-\varepsilon,\varepsilon)$,
where $\widetilde{M}_1$ is an integral submanifold for $T_\lambda$
through $p_0$ and $I$ corresponds to the integral curves of the unit
vector field $Y_1\in T_\mu$ on $U$. Moreover $\widetilde{M}_1$ is a
totally geodesic submanifold in $U$ and the integral curves of $Y_1$
are geodesics in $U$. We note that $Y_1$ is a parallel vector field
on $U$.

In the following, we shall prove that the integral curves of $Y_1$,
thought of as curves in $\mathbb{R}^{n+1}$, are circles of radius
$1/\sqrt 2$, all lying in parallel $2$-planes. In order to prove
this, consider $\{H/|H|,\eta_a\}_{a=1}^{k-1}$ to be an orthonormal
frame field in the normal bundle and $\{X_\alpha\}_{\alpha=1}^{m-1}$
an orthonormal frame field in $T_\lambda$, on $U$. We have
\begin{eqnarray*}
\trace B(A_H(\cdot),\cdot)&=&\sum_{\alpha=1}^{m-1} B(A_H(X_\alpha),X_\alpha)+B(A_H(Y_1),Y_1),\\
&=&\lambda mH-2\lambda B(Y_1,Y_1).
\end{eqnarray*}
This, together with \eqref{eq: caract_bih_Hparallel_I} and
 \eqref{eq: curv_Okumura}, leads to
\begin{equation}\label{eq: B(Y1,Y1)}
B(Y_1,Y_1)=-\frac{1}{\lambda}H,
\end{equation}
so $|B(Y_1,Y_1)|=1$. From here, since $A_{\eta_a}$ and $A_H$
commute, we obtain
\begin{equation}\label{eq: A_eta_Y1=0}
A_{\eta_a}(Y_1)=0,\qquad \forall a=1,\ldots,k-1.
\end{equation}
We also note that
\begin{equation}\label{eq: nabla B(Y1,Y1)}
\nabla_{Y_1}^{\mathbb{S}^n} B(Y_1,Y_1)=-\frac{1}{\lambda}(\nabla^\perp_{Y_1}H-A_H(Y_1))=-Y_1.
\end{equation}

Consider $c:I\to U$ to be an integral curve for $Y_1$ and denote by
$\gamma:I\to\mathbb{S}^n$, $\gamma=\imath\circ c$, where
$\imath:M\to\mathbb{S}^n$ is the inclusion map. Denote by
$E_1=\dot{\gamma}=Y_1\circ\gamma$. Since $Y_1$ is parallel, $c$ is a
geodesic on $M$ and, using equations \eqref{eq: B(Y1,Y1)} and
\eqref{eq: nabla B(Y1,Y1)}, we obtain the following Frenet equations
for the curve $\gamma$ in $\mathbb{S}^n$,
\begin{align}\label{eq: Frenet_gamma}
\nabla^{\mathbb{S}^n}_{\dot{\gamma}} E_1&=B(Y_1,Y_1)=-\frac{1}{\lambda}H=E_2,\nonumber\\
\nabla^{\mathbb{S}^n}_{\dot{\gamma}}E_2&=-E_1.
\end{align}

Let now $\widetilde{\gamma}=\jmath\circ\gamma:I\to \mathbb{R}^{n+1}$,
where $\jmath:\mathbb{S}^n\to\mathbb{R}^{n+1}$ denotes the inclusion map.
Denote by
$\widetilde{E}_1=\dot{\widetilde{\gamma}}=Y_1\circ\widetilde{\gamma}$.
From \eqref{eq: Frenet_gamma} we obtain the Frenet equations for
$\widetilde{\gamma}$ in $\mathbb{R}^{n+1}$,
\begin{align*}
\nabla^{\mathbb{R}^{n+1}}_{\dot{\widetilde{\gamma}}} \widetilde{E}_1&
=-\frac{1}{\lambda}H-\widetilde{\gamma}=\sqrt 2\widetilde{E}_2,\\
\nabla^{\mathbb{R}^{n+1}}_{\dot{\widetilde{\gamma}}}\widetilde{E}_2&=-\sqrt 2 \widetilde{E}_1,
\end{align*}
thus $\widetilde{\gamma}$ is a circle of radius $1/\sqrt 2$ in
$\mathbb{R}^{n+1}$ and it lies in a $2$-plane with corresponding
vector space generated by $\widetilde{E}_1(0)$ and
$\widetilde{E}_2(0)$.

Since $Y_1$ and $-\frac{1}{\lambda}H-\mathrm{x}$, with $\mathrm{x}$
the position vector field, are parallel in $\mathbb{R}^{n+1}$ along
any curve of $\widetilde{M}_1$, we conclude that the $2$-planes
determined by the integral curves of $Y_1$ have the same
corresponding vector space, thus are parallel.

Consider the immersions
$$
\phi: \widetilde{M}_1\times I\rightarrow\mathbb{S}^n,
$$
and
$$
\widetilde{\phi}=\jmath\circ \phi: \widetilde{M}_1\times I\to\mathbb{R}^{n+1}.
$$
Using the fact that $\widetilde{M}_1$ is an integral submanifold of
$T_\lambda$ and \eqref{eq: A_eta_Y1=0}, it is not difficult to
verify that $\widetilde{B}(X,Y)=0$, for all $X\in
C(T\widetilde{M}_1)$ and $Y\in C(TI)$, thus we can apply Lemma
\ref{lem: Moore lemma}. As the $2$-planes determined by the integral
curves of $Y_1$ have the same corresponding vector space and by
Corollary \ref{cor: non-full}, we obtain the orthogonal
decomposition
\begin{equation}\label{eq: decomp_Rn+1}
\mathbb{R}^{n+1}=\mathbb{R}^{n-1}\oplus \mathbb{R}^2
\end{equation}
and $U=M_1\times M_2$, where $M^{m-1}_1\subset \mathbb{R}^{n-1}$ and
$M_2\subset \mathbb{R}^2$ is a circle of radius $1/\sqrt 2$. We can
see that the center of this circle is the origin of $\mathbb{R}^2$.
Thus $M_1\subset\mathbb{S}^{n-2}(1/\sqrt 2)\subset \mathbb{R}^{n-1}$
and from Theorem \ref{th:hipertor}, since $U$ is biharmonic in
$\mathbb{S}^n$, we conclude that $M_1$ is a minimal submanifold in
$\mathbb{S}^{n-2}(1/\sqrt 2)\subset \mathbb{R}^{n-1}$. Consequently, the announced result holds locally.

In order to prove the global result we use the connectedness of $M$. Let $p\in M$ and let $V$ be an open neighborhood of $p$ given by the de Rham Theorem, as above, such that $U\cap V\neq \emptyset$. Consider $c_U$ and $c_V$ two integral curves for $T_\mu$, such that $c_U$ lies in $U$ and $c_V$ lies in $V$ and $c_U\cap c_V\neq\emptyset$. It is clear that the $2$-plane in $\mathbb{R}^{n+1}$ where $c_U$ lies coincides with the $2$-plane where $c_V$ lies. Therefore, the decomposition \eqref{eq: decomp_Rn+1} does not depend on the choice of $p_0$.

We can thus conclude that $M$ is an open part of a standard product
$$
M_1\times \mathbb{S}^{1}(1/\sqrt 2)\subset\mathbb{S}^{n},
$$
where $M_1$ is a minimal submanifold in $\mathbb{S}^{n-2}(1/\sqrt 2)$.

\end{proof}

By a standard argument, using the universal covering, we also obtain
the following result.

\begin{corollary}
Let $M^m$, $m>2$, be a proper biharmonic submanifold with
parallel mean curvature vector field in $\mathbb{S}^n$ and $|H|\in(0,1)$. Assume that $M$ is complete.
Then $|H|\in(0,\frac{m-2}{m}]$ and $|H|=\frac{m-2}{m}$ if and only if
$$
M=M_1^{m-1}\times\mathbb{S}^{1}(1/\sqrt 2)\subset\mathbb{S}^{n},
$$
where $M_1$ is a complete minimal submanifold of $\mathbb{S}^{n-2}(1/\sqrt 2)$.
\end{corollary}

If we consider the case of hypersurfaces, the condition on the mean
curvature vector field to be parallel is equivalent to the condition
on the mean curvature to be constant and Theorem \ref{th:
classif_parallelH} leads to the following result.

\begin{corollary}
Let $M^m$, $m>2$, be a proper biharmonic constant mean
curvature hypersurface with $|H|\in(0,1)$ in $\mathbb{S}^{m+1}$. Then $|H|\in(0,\frac{m-2}{m}]$. Moreover, $|H|=\frac{m-2}{m}$ if and only if $M$ is an open part of $\mathbb{S}^{m-1}(1/\sqrt{2})\times\mathbb{S}^1(1/\sqrt{2})$.
\end{corollary}

\begin{proof}
We recall that $|H|=\frac{m-2}{m}$ if and only if $\nabla A_H=0$ and
the principal curvatures of $M$ in the direction of $H$ are
constant, one of multiplicity $1$ and one of multiplicity $m-1$.
This implies that $M$ is an isoparametric hypersurface and, using a
result in \cite{BMO08, IIU08}, we conclude.
\end{proof}

\subsection{Biharmonic submanifolds with $\nabla^\perp H=0$ and $\nabla A_H=0$ in spheres}\quad

Inspired by the case $|H|=\frac{m-2}{m}$ of Theorem \ref{th: classif_parallelH}, in the following we shall study proper biharmonic submanifolds in $\mathbb{S}^n$ with parallel mean curvature vector field and parallel Weingarten operator associated to the mean curvature vector field.

We shall also need the following general result.
\begin{proposition}\label{prop: 2_princ_curv}
Let $M^m$ be a submanifold in $\mathbb{S}^n$ with nowhere zero mean
curvature vector field. If $\nabla^\perp H=0$, $\nabla A_H=0$ and
$A_H$ is orthogonal to $A_{\eta}$, for all $\eta\in C(NM)$,
$\eta\perp H$, then $M$ is either pseudo-umbilical, or it has two
distinct principal curvatures in the direction of $H$. Moreover, the
principal curvatures in the direction of $H$ are solutions of the
equation
\begin{equation}\label{eq: 2_princ_curv}
m t^2+\left(m-\frac{|A_H|^2}{|H|^2}\right)t-m|H|^2=0.
\end{equation}
\end{proposition}

\begin{proof} As $\nabla A_H=0$, the principal curvatures in the
direction of $H$ are constant. Denote by $\{X_i\}_{i=1}^m$ a local
orthonormal frame field on $M$ such that $A_H(X_i)=\lambda_iX_i$,
$i=1,\ldots,m$. Clearly, $\sum_{i=1}^m \lambda_i=m|H|^2$.

Since $A_H$ is parallel, $\nabla_X A_H(Y)=A_H(\nabla_X Y)$, thus
$R(X,Y)$ and $A_H$ commute for all $X, Y\in C(TM)$. In particular,
$$
R(X_i,X_j)A_H(X_j)= A_H(R(X_i,X_j)X_j),
$$
and by considering the scalar product with $X_j$ and using the
symmetry of $A_H$, we get
\begin{equation}\label{eq: aiajI}
(\lambda_i-\lambda_j)\langle R(X_i,X_j)X_j,X_i\rangle=0,\quad \forall i,j=1,\ldots,m.
\end{equation}
Consider $\{H/|H|,\eta_{a}\}_{a=1}^{k-1}$, $k=n-m$, a local
orthonormal frame field in the normal bundle of $M$ in
$\mathbb{S}^n$. We have
\begin{equation}\label{eq: B(Xi,Xi)}
B(X_i,X_i)=\frac{\lambda_i}{|H|^2}H+\sum_{a=1}^{k-1}\langle A_{\eta_a}(X_i),X_i\rangle\eta_a,
\end{equation}
and for $\lambda_i\neq \lambda_j$, as $X_i$ is orthogonal to $X_j$ and $A_H\circ
A_{\eta_a}=A_{\eta_a}\circ A_H$, for all $a=1,\ldots, k-1$, we
obtain
\begin{equation}\label{eq: B(Xi,Xj)}
B(X_i,X_j)=\frac{1}{|H|^2}\langle A_H(X_i),X_j\rangle H
+\sum_{a=1}^{k-1}\langle A_{\eta_a}(X_i),X_j\rangle\eta_a=0.
\end{equation}
By using \eqref{eq: B(Xi,Xi)} and \eqref{eq: B(Xi,Xj)} in the Gauss
equation for $M$ in $\mathbb{S}^n$, one gets
\begin{equation}\label{eq: aiajII}
\langle R(X_i,X_j)X_j,X_i\rangle=1+\frac{\lambda_i\lambda_j}{|H|^2}+
\sum_{a=1}^{k-1}\langle A_{\eta_a}(X_i),X_i\rangle\langle A_{\eta_a}(X_j),X_j\rangle.
\end{equation}
In fact, \eqref{eq: aiajI}, together with \eqref{eq: aiajII}, implies
\begin{equation}\label{eq: aiaj}
(\lambda_i-\lambda_j)\left(1+\frac{\lambda_i\lambda_j}{|H|^2}+\sum_{a=1}^{k-1}\langle
A_{\eta_a}(X_i),X_i\rangle\langle A_{\eta_a}(X_j),X_j\rangle\right)=0, \forall i,j=1,\ldots, m.
\end{equation}
Summing on $i$ in \eqref{eq: aiaj} we obtain
\begin{eqnarray*}
0&=&m|H|^2-\left(m-\frac{|A_H|^2}{|H|^2}\right)\lambda_j-m\lambda_j^2+\sum_{a=1}^{k-1}\langle
A_{\eta_a},A_H\rangle\langle A_{\eta_a}(X_j),X_j\rangle\\
&&-\sum_{a=1}^{k-1}\trace A_{\eta_a}\langle
A_{\eta_a}(X_j),A_H(X_j)\rangle.
\end{eqnarray*}
Since $\trace A_{\eta_a}=0$ and $\langle A_H, A_{\eta_a}\rangle=0$,
for all $a=1,\ldots, k-1$, we conclude the proof.

\end{proof}

\begin{corollary}\label{cor: bih_2curv}
Let $M^m$, $m>2$, be a proper biharmonic submanifold in
$\mathbb{S}^n$. If $\nabla^\perp H=0$, $\nabla A_H=0$ and $|H|\in
(0,\frac{m-2}{m})$, then $M$ has two distinct principal curvatures
$\lambda$ and $\mu$ in the direction of $H$, of different
multiplicities $m_1$ and $m_2$, respectively, and
\begin{equation}\label{eq: values lambda mu}
\left\{
  \begin{array}{ll}
    \lambda=\frac{m_1-m_2}{m}\\
    \mu=-\frac{m_1-m_2}{m}\\
    |H|=\frac{|m_1-m_2|}{m}.
  \end{array}
\right.
\end{equation}
\end{corollary}

\begin{proof}
Since $M$ is proper biharmonic, all the hypotheses of Proposition
\ref{prop: 2_princ_curv} are satisfied. Taking into account
\eqref{eq: caract_bih_Hparallel_II}, from \eqref{eq: 2_princ_curv}
follows that the principal curvatures of $M$ in the direction of $H$ satisfy the equation $t^2=|H|^2$. As $|H|\in (0,\frac{m-2}{m})$, $M$ cannot be pseudo-umbilical, thus it has two distinct principal curvatures $\lambda=-\mu\neq 0$ in the direction of $H$. If $m_1$ denotes the multiplicity of $\lambda$ and $m_2$ the multiplicity of $\mu$, from $\trace A_H=m|H|^2$ we have $(m_1-m_2)\lambda=m\lambda^2$. Since $\lambda\neq 0$, we obtain \eqref{eq: values lambda mu}. Notice also that $m_1\neq m_2$.
\end{proof}

We are now able to prove the following result.
\begin{theorem}\label{th: classif_parallel_H_AH}
Let $M^m$, $m>2$, be a proper biharmonic submanifold in
$\mathbb{S}^n$ with $\nabla^\perp H=0$, $\nabla A_H=0$ and $|H|\in
(0,\frac{m-2}{m})$. Then, locally,
$$
M=M^{m_1}_1\times M^{m_2}_2\subset \mathbb{S}^{n_1}
(1/\sqrt 2)\times \mathbb{S}^{n_2}(1/\sqrt 2)\subset\mathbb{S}^n,
$$
where $M_i$ is a minimal submanifold of $\mathbb{S}^{n_i}(1/\sqrt
2)$, $i=1,2$, $m_1+m_2=m$, $n_1+n_2=n-1$.
\end{theorem}

\begin{proof}
We are in the hypotheses of Corollary \ref{cor: bih_2curv}, thus
$A_H$ has two distinct eigenvalues $\lambda=\frac{m_1-m_2}{m}$ and
$\mu=-\frac{m_1-m_2}{m}$. Consider the distributions
\begin{align*}
T_\lambda&=\{X\in TM: A_H(X)=\lambda X\},\qquad \dim T_\lambda=m_1,\\
T_\mu&=\{X\in TM: A_H(X)=\mu X\},\qquad \dim T_\mu=m_2.
\end{align*}
As $A_H$ is parallel, $T_\lambda$ and $T_\mu$ are mutually
orthogonal, smooth, involutive and parallel, and from the de Rham
decomposition theorem follows that for every $p_0\in M$ there exists
a neighborhood $U\subset M$ which is isometric to a product
$\widetilde{M}_1^{m_1}\times \widetilde{M}_2^{m_2}$, such that the
submanifolds which are parallel to $\widetilde{M}_1$ in
$\widetilde{M}_1\times\widetilde{M}_2$ correspond to integral
submanifolds for $T_\lambda$ and the submanifolds which are parallel
to $\widetilde{M}_2$ correspond to integral submanifolds for
$T_\mu$.

Consider the immersions
$$
\phi: \widetilde{M}_1\times \widetilde{M}_2\rightarrow\mathbb{S}^n,
$$
and
$$
\widetilde{\phi}=\jmath\circ \phi: \widetilde{M}_1\times \widetilde{M}_2\to\mathbb{R}^{n+1}.
$$
It can be easily verified that $\widetilde{B}(X,Y)=B(X,Y)$, for all
$X\in C(T\widetilde{M}_1)$ and $Y\in C(T\widetilde{M}_2)$. Since
$A_H\circ A_\eta=A_\eta\circ A_H$ for all $\eta\in C(NM)$, we have that $T_\lambda$ and $T_\mu$ are invariant subspaces for $A_\eta$, for all $\eta\in C(NM)$, thus
$$
\langle B(X,Y),\eta\rangle=\langle A_\eta(X),Y\rangle=0,
\quad\forall \eta\in C(NM).
$$
Thus $\widetilde{B}(X,Y)=0$, for all $X\in C(T\widetilde{M}_1)$ and
$Y\in C(T\widetilde{M}_2)$, and we can apply Lemma \ref{lem: Moore
lemma}. In this way we have an orthogonal decomposition
$\mathbb{R}^{n+1}
=\mathbb{R}^{n_0}\oplus\mathbb{R}^{n_1+1}\oplus\mathbb{R}^{n_2+1}$ and
$\widetilde{\phi}$ is a product immersion. From Corollary \ref{cor:
non-full}, since $|H|\neq 1$, follows that $n_0=0$. Thus
$$
\widetilde{\phi}=\widetilde{\phi}_1\times\widetilde{\phi}_2:
\widetilde{M}_1\times\widetilde{M}_2\to \mathbb{R}^{n_1+1}\oplus\mathbb{R}^{n_2+1}.
$$
We denote by $M_1=\widetilde{\phi}_1(\widetilde{M}_1)\subset
\mathbb{R}^{n_1+1}$,
$M_2=\widetilde{\phi}_2(\widetilde{M}_2)\subset\mathbb{R}^{n_2+1}$ and
we have $U=M_1\times M_2\subset \mathbb{S}^n$.

Consider now $\{X_\alpha\}_{\alpha=1}^{m_1}$ an orthonormal frame
field in $T_\lambda$ and $\{Y_\ell\}_{\ell=1}^{m_2}$ an orthonormal
frame field in $T_\mu$, on $U$. From \eqref{eq:
caract_bih_Hparallel_I}, by using the fact that
$\lambda=-\mu=\frac{m_1-m_2}{m}$, we obtain
\begin{equation}\label{eq: H1H2}
\sum_{\alpha=1}^{m_1}B(X_\alpha,X_\alpha)=\frac{m_1}{\lambda}H,
\qquad \sum_{\ell=1}^{m_2} B(Y_\ell,Y_\ell)=-\frac{m_2}{\lambda}H.
\end{equation}
Since $\nabla^\perp H=0$, from \eqref{eq: H1H2} follows that
$M_1\times\{p_2\}$ is pseudo-umbilical and with parallel mean
curvature vector field in $\mathbb{R}^{n+1}$, for any $p_2\in M_2$.
But $M_1\times\{p_2\}$ is included in
$\mathbb{R}^{n_1+1}\times\{p_2\}$ which is totally geodesic in
$\mathbb{R}^{n+1}$, thus $M_1$ is pseudo-umbilical and with parallel mean curvature vector field in $\mathbb{R}^{n_1+1}$. This implies that $M_1$ is minimal in $\mathbb{R}^{n_1+1}$ or minimal in a hypersphere of $\mathbb{R}^{n_1+1}$. The first case leads to a contradiction, since $M_1\times\{p_2\}\subset \mathbb{S}^n$ and cannot be minimal in $\mathbb{R}^{n+1}$. Thus $M_1$ is minimal in a hypersphere $\mathbb{S}^{n_1}_{c_1}(r_1)\subset\mathbb{R}^{n_1+1}$, where $c_1\in\mathbb{R}^{n_1+1}$ denotes the center of the hypersphere.

In the following we will show that $c_1=0$. Since $U\subset \mathbb{S}^n$ and $M_1\subset\mathbb{S}^{n_1}_{c_1}(r_1)$, we get
$|p_1|^2+|p_2|^2=1$ and $|p_1-c_1|^2=r_1^2$, for all $p_1\in M_1$. This implies $\langle p_1,c_1\rangle=$ constant for all $p_1\in M_1$. Thus $\langle u_1,c_1\rangle=0$, for all $u_1\in T_{p_1}M_1$ and for all $p_1\in M_1$. From Lemma \ref{lem: Moore lemma} follows that $c_1=0$, thus $M_1\subset\mathbb{S}^{n_1}(r_1)\subset\mathbb{R}^{n_1+1}$.

From \eqref{eq: H1H2} also follows that the mean curvature of
$M_1\times\{p_2\}$ in $\mathbb{S}^n$ is $1$, so its mean curvature
in $\mathbb{R}^{n+1}$ is $\sqrt 2$. As
$\mathbb{R}^{n_1+1}\times\{p_2\}$ is totally geodesic in
$\mathbb{R}^{n+1}$ it follows that the mean curvature of $M_1$ in
$\mathbb{R}^{n_1+1}$ is $\sqrt 2$ too. Further, as $M_1$ is minimal in $\mathbb{S}^{n_1}(r_1)$, we get $r_1=1/\sqrt 2$.

Analogously, $M_2$ is minimal in a hypersphere $\mathbb{S}^{n_2}(1/\sqrt 2)$ in $\mathbb{R}^{n_2+1}$, and we
conclude the proof.

\end{proof}

\begin{corollary}
Let $M^m$, $m>2$, be a complete proper biharmonic submanifold in
$\mathbb{S}^n$ with $\nabla^\perp H=0$, $\nabla A_H=0$ and $|H|\in
(0,\frac{m-2}{m})$. Then,
$$
M=M^{m_1}_1\times M^{m_2}_2\subset \mathbb{S}^{n_1}(1/\sqrt 2)
\times \mathbb{S}^{n_2}(1/\sqrt 2)\subset\mathbb{S}^n,
$$
where $M_i$ is a complete minimal submanifold of $\mathbb{S}^{n_i}(1/\sqrt 2)$, $i=1,2$, $m_1+m_2=m$, $n_1+n_2=n-1$.
\end{corollary}

\begin{remark}
In the case of a non-minimal hypersurface the hypotheses $\nabla^\perp H=0$ and $\nabla A_H=0$ are equivalent to $\nabla^\perp B=0$, i.e. the hypersurface is parallel. Such hypersurfaces have at most two principal curvatures and the proper biharmonic hypersurfaces with at most two pricipal curvatures in $\mathbb{S}^n$ are those given by \eqref{eq: small_hypersphere} and \eqref{eq: product_spheres} (see \cite{BMO08}).
\end{remark}

\end{document}